\documentclass[12pt]{article}
\usepackage{amsfonts,amsthm}

\newcommand{\cent}[2]{\mathop{{\rm center}_{#1}}#2}
\newcommand{\codim}{\mathop{\rm codim}}
\newcommand{\Supp}[1]{\mathop{\rm Supp}#1}

\title{LETTERS OF A BI-RATIONALIST\\
{IV.~Geometry of log flips}}

\date{October 13th, 2000/January 2nd, 2001}

\author{V.V. Shokurov
\thanks{Partially supported by NSF
grants DMS-9800807 and DMS-0100991.}}

\theoremstyle{definition}
\newtheorem{adv}{Advertisement}
\newtheorem{conjecture}{Conjecture}
\newtheorem*{definition}{Definition}

\theoremstyle{remark}
\newtheorem{example}{Example}
\newtheorem{remark}{Remark}
\newtheorem{question}{Question}
\newtheorem{warning}{Warning}

\theoremstyle{theorem}
\newtheorem{corollary}{Corollary}
\newtheorem*{lemma}{Lemma}
\newtheorem*{monot}{Monotonicity}
\newtheorem*{theorem}{Theorem}

\begin{document}

\maketitle

Flops and flips first appeared
in mathematics as geometrical
constructions:
\begin{description}
\item{(1)}
during Fano's modification of
a $3$-fold cubic into
a Fano $3$-fold $X_{14}\subset
\mathbb P^9$ \cite[Theorem 4.6.6]{IP};
\item{(2)}
Atiyah's flop: one of his first papers \cite{at} in 1958
treated the simultaneous resolution
of the surface ODP, and was the initial stimulus for Brieskorn's
simultaneous resolution of Du Val singularities\footnote{Thanks
to Miles Reid for this historical remark. He added also
{``}Possibly a little later,
Moishezon (and Hironaka) were using the same
kind of thing to construct algebraic spaces (minischemes) that were not
varieties.

However, as I said in my Old Person's View, one can trace the idea back
through Zariski and Kantor and Cremona, even as far back as papers of
Beltrami in 1863 and Magnus in 1837 referred to in Hilda Hudson's
bibliography -- these papers study the standard monoidal involution
$\mathbb P^3 -\to \mathbb P^3$
given by $(x,y,z,t) \mapsto (1/x, 1/y, 1/z, 1/t)$, which
flops the 6 edges of the coordinate tetrahedron".};
\item{(3)}
Kulikov's perestroikas
\cite[Modifications in 4.2-3]{Ku};
\item{(4)}
Francia's flip (see Example~\ref{frc} below);
\item{(5)}
Reid's pagodas \cite{R};
\item{(6)}
semistable flips \cite{T}
\cite{K88} \cite{Sh94};
\item{(7)}
Kawamata's nonsingular $4$-fold flip
\cite{K89};
\item{(8)}
geometrical $4$-fold flips \cite{Ka};
and
\item{(9)}
the Thaddeus principle \cite{Th}
\cite{BH}.
\end{description}
However in general, for higher
dimensions, one can hardly imagine
an effective and explicit
{\em geometric\/} construction
(for instance, a chain of certain
blow-ups and blow-downs)
for {\em flips\/}, even for log ones,
except for very special situations with
extra structures, e.g., as in (6) (8), and for
moduli spaces as in (9).
On the other hand, we hope that
the log flips exist and this can be
established in a more formal and algebraic
way.
Recently, this was done for
the log flips up to dimension $4$
\cite[Corollary~1.8]{Sh00}.
Since these flips were obtained without the use
of any classification or
concrete geometry of them,
it is worthwhile in the
aftermath to get some of
the aforementioned geometrical facts.
This is a goal of the note which
we pursue in
a more general situation.
For the convenience of the reader, we
put the list of notation and terminology
at the end of the letter.

\begin{definition}
A birational transform $X-\to X^+/Z$ between
\label{qflip}
two birational contractions $f:X\to Z$ and
$f^+:X^+\to Z$ of
normal algebraic varieties is called
a (directed) $D$-{\em quasi-flip}$/Z$ or, shortly,
-{\em qflip}$/Z$, for
a Weil $\mathbb R$-divisor $D$ on $X$,
when there is a semiample$/Z$
$\mathbb R$-Cartier divisor $D^+$ on $X^+$ such
that
$f_*^+D^+\sim_{\mathbb R} f_*D$.

A $D$-qflip$/Z$ can be given in a {\em log form} or, shortly,
in  {\em lf}, that is,
in terms of log structures on $X/Z$ and
$X^+/Z$, namely, there are
Weil $\mathbb R$-divisors $B$ and
$B^+$ on $X$ and $X^+$ respectively such that:
\begin{itemize}
\item
$f_*^+B^+=f_*B$; and
\item
$D=K+B$ and $D^+=K_{X^+}+B^+$.
\end{itemize}
The qflip is a {\em log} one if in addition:
\begin{itemize}
\item
$B$ and $B^+$ are boundaries; and
\item
pairs $(X,B)$ and $(X^+,B^+)$ are log canonical.
\end{itemize}
\end{definition}

Note that up to an $\mathbb R$-linear equivalence of
$D$ and/or $D^+$ we can assume that
$f_*^+D^+= f_*D$ in the definition.
Then any $D$-qflip is a qflip in lf for
some $B$ and $B^+$ (but maybe not a log qflip).
(We always take all canonical divisors $K, K_{X^+}$, etc.
on modifications of $X$ given by the same differential form,
or by the same bi-divisor \cite[Example~1.1.3]{Sh96}.)
Any $D$-flip is a $D$-qflip with
$D^+$ as the birational transform of $D$.
The inverse holds when $X^+/Z$ is small and
$D^+$ is ample$/Z$; for example,
by Monotonicity below, for any qflip in lf, $X^+/Z$ is
small if $X/Z$ is small, $-D=-(K+B)$ is nef$/Z$, and
$(X,B)$ is terminal in codimension $2$, that is,
$a(Y)>1$ for any subvariety $Y\subset X$
of codimension $\ge 2$ in notation below.
Thus a log qflip is a log flip under the last
conditions, and with ample $D^+=K_{X^+}+B^+/Z$.
Log qflips, with nonsmall $X^+/Z$ and
a boundary $B^+$, are naturally
induced by log flips on the reduced part
of $B$ (cf. the proof of \cite[Special
termination~2.2]{Sh00}).

Even always assuming
that the characteristic
of base field $k$ is $0$, we
expect that most of the results and
statements below hold without such
an assumption, e.g.,
the following generalizations of
\cite[Lemma~5-1-17]{KMM} and
Monotonicity \cite[(2.13.3)]{Sh83}
 -- our basic tools.

\begin{lemma}
Let $X-\to X^+/Z$ be a $D$-qflip
for an $\mathbb R$-Cartier divisor $D$ on $X$
such that:
\begin{enumerate}
\item
$X/Z$ is a $D$-{\rm contraction, that is,
$-D$ is numerically ample$/Z$\/}
{\rm \cite[Section~5]{Sh96}\/}, and
\item
$X/Z$ is a nonisomorphisms.
\end{enumerate}
Then
$$
{^+}c\le d+1
$$
where
\begin{itemize}
\item
$d=d(X/Z)$ is the {\rm minimal\/} dimension
of the irreducible components of
the exceptional locus $E$ of $X/Z$; and
\item
${^+}c$ is the {\rm minimal\/} codimension
in $X^+$ of the irreducible components of
the {\rm rational transform\/}
${^+}E$ of $E$ in $X^+/Z$.
\end{itemize}
If $E$ has the {\rm pure\/} dimension $d${\rm , that is, 
each irreducible component of $E$ has
the dimension $d$}, then ${^+}c$ can be taken
as the {\rm maximal\/} codimension.
\end{lemma}

\begin{warning}
\label{incl}
In general, ${^+}E$ is quite different from
the exceptional locus $E^+=E(X^+/Z)$.
However $E^+\subseteq {^+}E$ whenever
$X-\to X^+$ is an isomorphism on $X\setminus E$;
for instance, the latter holds for the $D$-flips
of $D$-contractions but not for all qflips.
\end{warning}

\begin{remark}
\label{max}
The minimal dimensions and codimensions
can be replaced in the dual form of
the lemma by maximal ones, namely,
$c\le {^+}d+1$ in its {\em maximal\/} form.
The lemma itself, in the maximal form, is dual to
its symmetric statement $c\le {^+}d+1$ in
the {\em minimal\/} form holding
for the same $E$ and ${^+}E$.

Note also that taking the minimal dimension and
maximal codimension we consider only
nonempty components, in particular, such
(co-)dimensions are
defined only for nonempty subvarieties.
This explains Condition 2.
\end{remark}

\begin{proof}
After a birational contraction of $X^+/Z$ given
by $D^+$ we can assume that $X^+/Z$ is
$-D^+$-contraction; this change only
increases ${^+}c$.

Then the ampleness of $p_1^*(-D)+p_2^*D^+$ on
$X\times_Z X^+/Z$ and \cite[Proof of Lemma~5-1-17]{KMM}
imply that $X\times_Z X^+$ is divisorial
over $Z$ (see also \cite[Negativity 1.1]{Sh92}).
Moreover, for each irreducible component $Y$ of $E$ and
its rational image ${^+}Y\subseteq {^+}E\subset X^+$,
$$
\dim Y+\dim {^+}Y\ge \dim Y\times_Z {^+}Y=\dim X-1.
$$
That gives the required inequality.

The last statement for the pure $d$ follows
from the maximal case mentioned in
Remark~\ref{max}.
\end{proof}

\begin{monot}
Let $(X,B)-\to (X^+,B^+)/Z$ be
a qflip in lf with nef $-(K+B)/Z$.
Then, for each prime bi-divisor $P$ of $X$,
$$
a(X^+,B^+,P)\ge a(X,B,P).
$$
Moreover,
$$
a(X^+,B^+,P)> a(X,B,P)
$$
for each $P$ with $\cent{X}{P}\subseteq E$
when $-(K+B)$ is numerically ample$/Z$;
the equivalent inclusion is
$\cent{X^+}{P}\subseteq {^+E}$.
\end{monot}

\begin{proof}
As for \cite[(2.13.3)]{Sh83}.
\end{proof}

Let $(X/Z,B)$ be a log pair
with a boundary $B$ such
that:
\begin{description}
\item{\rm (BIR)}
\label{bir}
$f:X\to Z$ is a {\em birational\/} contraction which
we always consider locally over some fixed
point in $Z$;
and
\item{\rm (WLF)}
\label{wlf}
the pair is a {\em weak log canonical Fano\/}
contraction, that is,
$(X,B)$ is log canonical, and
$-(K+B)$ is nef$/Z$; 
\end{description}
it is said to
be a {\em log canonical Fano} contraction or
{\em log contraction\/} when
$-(K+X)$ is numerically ample$/Z$.
Note that $-(K+B)$ is big$/Z$ by (BIR).
The log pairs include, in particular,
the birational {\em log contractions\/} of
LMMP (Log Minimal Model Program), which are
birational contractions $X/Z$ of
extremal faces numerically
negative$/Z$ with respect to $K+B$
\cite[5.1.1b]{Sh96}.
However in this letter we do not touch
fibred contractions \cite{AW}.

The most fundamental questions in
geometry concern {\em dimensions\/}.
In our situation they are
\begin{itemize}
\item
$n=\dim X$; and
\item
the {\em minimal\/} dimension $d=d(X/Z)$
for the exceptional locus $E$ of $X/Z$.
\end{itemize}
Other more modern numerical invariants:
\begin{itemize}
\item\label{logl}
the {\em (log) length\/}
$l=l(X/Z,B)$ of
$(X/Z,B)$, that is the minimal
$-(K+B.C)$ for generic curves $C$ in
the covering families of
{\em contracted\/} locus $E$ (this is
the exceptional locus whenever
$X/Z$ is birational, and
$E=X$ otherwise) of $X/Z$; and
\item
the {\em m.l.d.
(minimal log discrepancy)\/}
$a=a(E)=a(X,B,E)$ of $(X,B)$ {\em in\/} $E$,
that is the minimal
log discrepancy $a(X,B,P)$ at
prime bi-divisors $P$
having the center {\em in\/} $E$;
the latter means that
$\cent{X}{P}\subseteq E$.
\end{itemize}
It is known that the dimension $d$
depends on the length \cite[Theorem]{Sh96b}, and
on the singularities \cite[Th\'eor\`eme~0]{Be}.
Sometimes a more subtle interaction
occurs.

\begin{example}
\label{benv}
Suppose that $(X,B)$ has only canonical
(terminal) singularities
in codimension $2$, and a curve $C$
is an irreducible component of
the exceptional locus $E$ of projective $X/Z$.
Then the {\em existence of log flips
in dimension $n\ge 3$ in the formal$/\overline{k}$, or
analytic category when $\overline{k}=\mathbb C$, implies
that $(K+B.C)\ge -1$ (respectively,
$>-1$)\/}.
More precisely, for $n\ge 2$, we can assume
just $a(C)\ge 1$ (respectively $>1$).
In other words,
this means that the length $l$ of
the contraction is $\le 1$
(respectively $< 1$)
whenever $d=1$ and $a(C)\ge 1$
(respectively $>1$).
One can drop the existence of
the log flip in
dimension $n\le 4$.

We verify that $l\le 1$
in the canonical case;
the terminal case is similar
(cf. \cite[Lemma~3.4]{ChP}).
Indeed, suppose that $l>1$ and
$a(C)\ge 1$.
Then over a small neighborhood of
$f(C)$ in the classical
complex topology for $\overline{k}=\mathbb C$
(or formally over arbitrary algebraically closed
$\overline{k}$), there exists
a rather generic hyperplane section $H/Z$ that
intersects $C$ transversely
in a single point.
In addition, changing the contraction
over such a neighborhood we can
assume that $E=C$.
So, locally$/Z$, $(X/Z,B+H)$ is again
a log pair under (BIR) and (WLF).
Since the exceptional locus $E=C$ is
a {\em proper\/} subvariety in $X$, the new m.l.d.
$a:=a(X,B+H,E)=a(X,B,E)\ge 1$, too.

On the other hand,
according to our assumptions,
there exists a flip $X-\to X^+/Z$
with respect to $K+B+H$.
Actually it is also the flip for
$K+B$ and is the $(-H)$-flip
\cite[Corollary~3.4]{Sh00}.
So, the flip transform ${^+H}$ of $H$ is
the birational transform $H^+$ of $H$ and
numerically negative$/Z$ on the exceptional locus
$E^+$ of the flipped contraction $X^+/Z$.
Hence $E^+\subseteq {^+}E\subset H^+$, and
$E^+={^+}E$ unless $n=2$ with $E^+=\emptyset$
(cf. Warning~\ref{incl}).
Moreover, ${^+}E$ has the minimal codimension
$c^+=2$ by the lemma, and
the m.l.d. ${^+}a=a(X^+,B^++H^+,{^+}E)\le 1$;
it is enough to establish the latter for $n=2$,
when it is well-known
\cite[Example~4.2.1]{Sh96}.
But this contradicts to
the assumption $a\ge 1$ because
$a<{^+}a\le 1$ by Monotonicity.
\end{example}

\begin{remark}
\label{cf}
In the last paragraph we proved
a little bit more.
Let $(X/Z,B)$ be a purely log terminal
pair with the reduced divisor $H$.
Then each flip of $(X/Z,B)$
with $a\ge 1$ gives
the flip on $(H/f(H),B_H)$ where
$B_H$ is given by the adjunction.
(So, then $d\ge 2$ by the lemma
when $H\sim_{\mathbb R}-h(K+B-H)$, with $h\in\mathbb R$,
is numerically ample and $\not\equiv 0/Z$
as in the example.)
Therefore, for a purely log terminal
and canonical in codimension $2$ pair
$(X/T,B)$,
LMMP with only {\em flipping\/}
contractions $(X/Z/T,B)$ induces LMMP
on $(H/f(H),B_H)$
(cf. the proof of \cite[Special Termination
2.3]{Sh00}).
Moreover, the same holds for any chain of
{\em birational\/} contractions in LMMP for
$(X/T,B)$ unless one of them contracts
a component of $\Supp{B}$.
\end{remark}

\begin{adv}
A generalization and applications of
the improvement in Remark~\ref{cf} will be treated in
one of the following letters.
\end{adv}

In the $3$-dimensional terminal case
with $B=0$, Example~\ref{benv} implies
the Benveniste result \cite[Th\'eor\`eme~0]{Be}\footnote{In
general the strict inequality in the theorem fails
in presence of canonical singularities along 
curve $C$, e.g., when $C$ is obtained by the contraction
along the second factor of a surface 
$\mathbb P^1\times \mathbb P^1$ in a nonsingular 3-fold $X$
with the normal bundle $\pi_1^*\mathcal O_{\mathbb P^1}(-1)
\oplus\pi_2^*\mathcal O_{\mathbb P^1}(-2)$ where
$\pi_i :\mathbb P^1\times \mathbb P^1\to \mathbb P^1$ 
is the projection on $i$-th factor.}; now
without linear systems arguments.
But it looks difficult to apply his
arguments in higher dimensions;
even in dimension $4$.
Deformation arguments in any
dimension gives the weaker inequality
$l<2$ in Example~\ref{benv}.
In general, $d>l/2$
even for more difficult singularities
\cite[Theorem]{Sh96b}.
On the other hand, we expect

\begin{conjecture}
\label{mineq}
Under conditions (BIR) and (WLF),
suppose that $X/Z$ is {\em projective\/}.
Then $d\ge a-1$
($>$ in the {\em log Fano\/} case), or,
equivalently, $d\ge\lceil a-1\rceil$.

In addition, if $d=\lceil a-1\rceil,
E\not=\emptyset$, and $(X/Z,B)$ is
a {\em log contraction, that is, $X/Z$ is a $D$-contraction
for $D=K+B$\/}, then,
for any log qflip $X-\to X^+/Z$,
the transform ${^+}E$
satisfies the following properties:
\begin{description}
\item{(CDM)}\label{cdm}
${^+}c= d+1$;
\item{(NSN)}\label{nsn}
each irreducible component of ${^+}E$ of
the minimal codimension $d+1$ is {\em nonsingular\/} as
a scheme point of $X^+$; and
\item{(PDM)}\label{pdm}
if $d$ is the {\em pure dimension\/}, then ${^+}E$ is
also of the pure codimension $d+1$.
\end{description}
Moreover, for the log flip $X-\to X^+/Z$,
${^+}E=E^+$ is the exceptional locus of $X^+/Z$.

It is enough to establish the conjecture when
$E\not=\emptyset$.
Otherwise we put $d=-\infty$ and $a=-\infty$ as it
used to be.
\end{conjecture}

\begin{example}
\label{wish}
In particular,
if $X$ is nonsingular and $B=0$
then $a=n-d$ and
Conjecture~\ref{mineq} implies that
$d\ge (n-1)/2$ (respectively $d> (n-1)/2$.
This is Wi\'sniewski's inequality
\cite[Theorem 1.1, p.~147]{W} in
a {\em single\/} formula).
\end{example}

Perhaps there is a relation
between the length $l$ and
the m.l.d. $a$.
I am not sure.
But definitely, $l/2\ge a-1$ does not
hold always.

\begin{example}
\label{frc}
Francia's flip corresponds to
the contraction, which is obtained from the relative model
$Y$ after contraction of its plane $E$ in \cite[Section~2]{F},
and it has $a=3/2$ and $l=1/2$.
So, the above inequality $d>l/2$ is
less sharp than that of in Conjecture~\ref{mineq}.
\end{example}

Conjecture~\ref{mineq} can be derived from
the conjecture on existence of log flips \cite[Conjecture~5.1.2]{Sh96}
and another conjecture on the m.l.d.
\cite[Problem~5a]{Sh88}:

\begin{conjecture}
\label{mld}
For any (scheme) point $P\in X$
the m.l.d. $a(X,B,P)\le
\codim P$,
with $=$ holds only when $P$ is nonsingular
in $X$ and $B=0$ near $P$.
Taking hyperplane sections, it
is enough to prove that
$a(X,B,P)\le \dim X$ and the $=$ case for
closed points.

Moreover, the nonsingularity of $X$ still holds
when we replace $\codim P$ by $\codim P-\varepsilon$
with any $0\le \varepsilon<1$.
Equivalently, $\lceil a(X,B,P)\rceil\le \codim P$
with $=$ holds only when $P$ is nonsingular in $X$.
\end{conjecture}

\begin{warning}
In the conjecture we still assume that
$B$ is a boundary!
\end{warning}

Conjecture \ref{mld} was proven for $\codim P\le 3$
(after a $\mathbb Q$-factorialization,
follows from \cite[Corollary~3.3]{Sh96} with a final step by
Markushevich \cite[Theorem~0.1]{M}).
The stronger form with $\lceil\quad\rceil$ follows from
the weaker one, the covering trick and LMMP.
From LMMP, we need only the existence
of $\mathbb Q$-factorializations.
(It is expected that the next case with $a(X,0,P)=\codim P-1$
corresponds to the higher dimensional cDV singularities;
cf. \cite[Proposition~3.3]{am}.)

Since the m.l.d. measures the singularity, it
is natural to expect that it decreases under
the specialization that is stated in Ambro's
conjecture \cite[Conjecture~0.2]{am}.
It implies Conjecture~2, and is proved
up to dimension $3$ \cite[Theorem~0.1]{am} and
for toric varieties by \cite[Theorem~4.1]{am}.
The former gives again $\codim P\le 3$ and the latter
gives the toric case.

\begin{example}
\label{tor}
Conjecture~2 holds for toric varieties with invariant $B$.
\end{example}

\begin{theorem}
The {\rm existence of log flips in
dimension\/} $n$ and
{\rm Conjecture~\ref{mld}\/} in dimension $m$ implies
{\rm Conjecture~\ref{mineq}\/} in dimension $n$
for any $d\le m-1$.

Actually, log flips can be weaken to log qflips.
More precisely,
it is enough to have the existence of
log qflips for log contractions $(X/Z,B)$
instead of log flips.
\end{theorem}

\begin{proof}
According to the closing remark in
Conjecture~\ref{mineq},
we can assume that $E\not=\emptyset$.
In particular, $n\ge 2$.

We can suppose also that $a>0$.
Otherwise, $a=0$ and
$d\ge 1>-1=a-1=\lceil a-1\rceil$.

Since $X/Z$ is projective, after perturbation of $B$ we
can assume that $X/Z$ is a log contraction.
Note that taking quite a small perturbation,
which increases $B$ and decreases $a$, we
preserve $\lceil a-1\rceil$.
If the original $B$ gave the log contraction
we do not change $B$.

Now we can apply the lemma.
Let $X-\to X^+/Z$ be a log qflip of $(X/Z,B)$
with a boundary $B^+$ on $X^+$.
Then by the lemma
there is an irreducible component $Y\subseteq {^+}E$
such that $\codim{Y}\le d+1$.
Hence Conjecture~\ref{mld} implies
that $a(X^+,B^+,Y)\le d+1$, and Monotonicity
$$
a=a(X,B,E)\le a(X^+,B^+,Y)\le d+1
$$
gives the required inequality.

Under additional assumptions in
Conjecture~\ref{mineq}, if $\codim{Y}\le d$ then,
according to the same inequalities,
$a\le d$ and $\lceil a-1\rceil\le d-1<d$.
This proves (CDM) because, according
to these assumptions, $d=\lceil a-1\rceil$ or
$\lceil a\rceil=\lceil a(X^+,B^+,Y)\rceil=d+1$.
So, (NSN) follows from Conjecture~\ref{mld}.
The last statement in the lemma implies (PDM).
If $X-\to X^+/Z$ is the log flip, then
$E^+= {^+}E$ (cf. Warning~\ref{max}).
Indeed, for small $X/Z$, the inverse transform
is the anti-flip.
Otherwise $a\le 1,d=n-1=\lceil a-1\rceil\le 0$,
and $n\le 1$, which contradicts to our assumptions.
\end{proof}

\begin{corollary}\label{torc}
Conjecture~\ref{mineq} holds for the toric contractions
$X/Z$, with only canonical singularities and $B=0$, 
which are numerically negative with respect to $K$.
\end{corollary}

A more general case we consider in one of our
future letters (see Advertisement~\ref{acclt}).
It would be interesting to know whether the
combinatorics behind this statement were known.
In particular, how important is the projectivity
in this statement (cf. Question~\ref{nonp})?

\begin{proof}
Immediate by Example~\ref{tor} and the existence of
toric flips \cite[Theorem~0.2]{R83}. 
We do not need to perturb $B=0$ since $K$ itself is
numerically negative. 
\end{proof}

\begin{corollary}
\label{small}
The theorem holds without {\rm Conjecture~\ref{mld}\/}
for all $d\le m=2$.
\end{corollary}

\begin{proof}
Immediate by \cite[Corollary~3.3]{Sh96} and
\cite[Theorem~0.1]{M}.
\end{proof}

The main inequality $d\ge a-1$ in
Conjecture~\ref{mineq} can be established for
rather high dimensions $d$ without flips.

\begin{example}
\label{big}
For all $d\ge (n-1)/2$, Conjecture~\ref{mineq}
follows form Conjecture~\ref{mld} for
$m\le (n+1)/2$;
in particular, up to $n=6$ we can drop
Conjecture~\ref{mld}.
Indeed, let $Y$ be an irreducible
component of $E$
then $a\le\codim{Y}\le (n+1)/2\le d+1$.
Moreover, $d=\lceil a-1\rceil$ only if
$d=(n-1)/2$, the dimension is pure and
$E$ is nonsingular in each of its irreducible
components as a scheme point;
the additional statements (CDM), (NSN) and
(PDM) in Conjecture~\ref{mineq}
hold by the lemma and our hypothesis
(cf. the proof of the theorem).
Otherwise $d\ge n/2$ is integral, and
$\lceil a\rceil\le n/2=d<d+1$.
\end{example}

\begin{corollary}
\label{offl}
In dimensions $n\le 6$ one can
drop {\rm Conjecture~\ref{mld}\/} in
the {\rm theorem\/}.
\end{corollary}

\begin{proof}
In Example~\ref{big} it was proven
for all $d\ge (n-1)/2$ because
$m=(n+1)/2\le (6+1)/2=7/2$.
For $d\le 2<(n-1)/2\le (6-1)/2=5/2$, the corollary 
was proven in Corollary~\ref{small}.
\end{proof}

\begin{corollary}
In dimension $n\le 4$ one can
drop both conjectural hypotheses
in the {\rm theorem, namely,
the existence of log flip and Conjecture~\ref{mld}\/}.
\end{corollary}

\begin{proof}
Immediate by Corollary~\ref{offl} because
the log flips exists.
The latter was proven in \cite[Corollary~1.8]{Sh00}
when $(X,B)$ is Kawamata log terminal.
The other log flips also exist  due
to \cite[Special termination~2.3]{Sh00} and
the local log semiampleness (cf. \cite[Conjecture~2.6]{Sh96}, and
see the proof of \cite[Log Flip Theorem~6.13]{Sh96}).

Actually, by Example~\ref{big} it is
enough to consider the case with pure $d=1$.
Then $a\le 2$.
Indeed, otherwise after a strict
log terminal
resolution we can assume that
each reduced component $H$ in $B$ is
$\mathbb Q$-Cartier, and intersects properly
the curves $C$ of $E$.
Indeed, we can construct the log terminal resolution$/X$
using Kawamata log terminal flips and
\cite[Special termination~2.3]{Sh00} as in
\cite[Reduction~6.5]{Sh92} (cf. the proof
of \cite[Proposition~10.6]{Sh00}).
This is not an isomorphism only over a finite set of
points in $C$ because $X$ is nonsingular in
the generic points of $C$ by our assumption.
This is impossible by the Kawamata log terminal case
(cf. Example~\ref{benv}) because
$C$ is contractible at least in the formal or
analytic category, $a>2$ can only
be increased after the construction,
and now the flip in $C$ exists
(see \cite[Remark~1.12]{Sh00}).
\end{proof}

\begin{example}[Minimal contractions]
\label{minc}
Under conditions (BIR) and (WLF),
suppose that $X/Z$ is projective and dimension
$d=\lceil a-1\rceil$ is {\em pure\/} (or maximal).
Then a log pair $(X/Z,B)$ will be called a {\em minimal
(log) contraction\/} (respectively, when $B\not=0$).

In particular, a minimal log contraction with
$d=a-1$ is possible only when $(X/Z,B)$ is a $0$-log pair, 
that is, $K+B\equiv 0/Z$ in our situation.
This follows from a more subtle version of
$>$ in Monotonicity under $\not\equiv 0/Z$
by \cite[Negativity~1.1]{Sh92}, or
LMMP including the log termination and
our theorem.
If $X/Z$ is projective under hypotheses of
the theorem (cf. its proof)
there exists a nonidentical directed flop
$X-\to X^+/Z$ with the transform ${^+}E$ of pure codimension and
satisfying the properties (CDM), (NSN), and (PDM)
of Conjecture~\ref{mineq} (perturb $B$ by $D$ as
a negative to a polarization).
In particular, such a flop is unique when
$X$ is $\mathbb Q$-factorial and $X/Z$ is
{\em formally extremal\/},
that is, the {\em formal\/}
(in the formal or analytic category) relative Picard number
of $X/Z$ is $1$.
Moreover, for dimension
$d=(n-1)/2=a-1$, it is expected that
$(X/Z,B)$ is nonsingular in
the irreducible components of $E$ as scheme points,
$f(E)$ is a {\em closed point\/} (take a general hyperplane
section of $f(E)$ when $\dim f(E)\ge 1$), and
each directed flop (qflop) $X-\to X^+/Z$
should be {``}symmetric{"}, that is,
the exceptional $E^+={^+E}$ of
pure dimension $d=(n-1)/2$ with the nonsingular
irreducible components of $E^+$ as scheme points
(cf. Questions~\ref{nonp} and \ref{minf} below).
The same follows from LMMP for any
projective flop $(X^+/Z,B)$ of $(X/Z,B)$ as a
composition of directed ones.

For instance, if $X$ is nonsingular and $B=0$,
then it is expected that $d\ge (n-1)/2$
(cf. Example~\ref{wish}).
So, by Conjecture~\ref{mld} $(X/Z,0)$ is
minimal only when $d=a-1=(n-1)/2$
and the dimension is pure.
Again by Conjecture~\ref{mld}
it is expected that any directed or/and
projective flop $X^+/Z$ is nonsingular with
the same number of irreducible components
of $E^+$ as for $E$ (the number
of exceptional prime divisors over $E$ or
$E^+$ with the log discrepancy $a=(n+1)/2$)
whenever
$X^+/Z$ is {\em formally\/} $\mathbb Q$-{\em factorial\/}
(in the formal or analytic category;
cf. Question~\ref{minf}  below).
The latter should hold for any nontrivial
flop when $X/Z$ is formally extremal.
In this case the flop is unique when it exists and
will be called {\em minimal formally extremal\/}.
The LMMP implies that each of the directed and projective
flops to $X^+/Z$ is a composition of
(formally) extremal contractions
and flops; only such flops are enough when both $X$ and $X^+/Z$ are
(formally) $\mathbb Q$-factorial.

An elementary example with nonsingular $X$, $B=0$, and 
$E=\mathbb P^d$ belongs to the toric geometry.
Its invariant divisors are $(n+1)/2=d+1$ numerically
negative $D_i^-$ (intersecting $E$ up to the linear equivalence
by a negative to its hyperplane;
their intersection in $X$ is $E$ itself) and
$(n-1)/2+1=d+1$ numerically positive $D_j^+$ (intersecting
$E$ in hyperplane sections in a general position --
an anti-canonical divisor in total).
The construction of such a contraction $X/Z$ see
in \cite[Example~3.12.2(iii)]{Isk} with
$r=(n+1)/2=d+1$, and $a_1=\dots=a_r=1$.
This contraction and its flop are
{\em formally or analytically toric\/}
(cf. a conjecture after the proof of
\cite[Theorem~6.4]{sh95}) and
{\em semistable\/} (cf. (6)); in this situation
the latter means that
there exists a nonsingular hypersurface $H$ passing through
$E$ and $\equiv 0/Z$.
Thus they induce a flop on this hypersurface
as in Example~\ref{amc} below with $E/pt$.
Moreover, it is a symplectic one (cf. Question~\ref{amin}).
Each toric contraction $(X/Z,B)$ is 
algebraically (analytically or formally) log $i$-symplectic
for any $0\le i\le n=\dim X$
with $B$ as invariant divisor
(take a general linear combination of invariant $i$-forms
$\wedge^i(d z_j)/z_j$).
In our situation, we can take $3$-form ($3$-symplectic structure)
$$
\frac{d(f_1^-f_1^+)}{f_1^-f_1^+}\wedge
(\sum_{2\le i,j\le d+1}\frac{d(f_i^-)d(f_j^+)}{f_i^-f_j^+})
$$
where $f_i^-$ and $f_j^+$ are respectively local sections$/Z$ of
$\mathcal O_X(-D_i^-)$ and $\mathcal O_X(-D_j^+)$.
Then after a generic perturbation of products
$f_1^-f_1^+$ and $f_i^-f_j^+$ into $g_{11}$, with
$H=\{g_{11}=0\}$ vanishing on $\mathbb P^d$, and
$g_{ij}$ (the latter does not vanishing on $\mathbb P^d$
at all),
the $3$-form
$$
\frac{dg_{11}}{g_{11}}\wedge
(\sum_{2\le i,j\le d+1}\frac{d(f_i^-)d(f_j^+)}{g_{ij}})
$$
induces by its residue
a symplectic $2$-form on $H$.
Such flops will be called {\em induced toric\/}.

Note that, according to A. Borisov (a private
communication) the {\em top\/} m.l.d. in dimension $n$ for
the toric {\em isolated $\mathbb Q$-factorial
singularities\/} $P$ is $a(P)=n/2<(n+1)/2$. 
The latter is the mld for
the minimal contraction of nonsingular $X$ with $B=0$.
Actually, {\em for any toric isolated $\mathbb Q$-Gorenstein
singularity $P$, $a(P)\le (n+1)/2$, and $a(P)=(n+1)/2$ only
for the above toric contraction.\/} 
For such a toric singularity $P\in Z$,  
a toric projective $\mathbb Q$-factorialization $X/Z$ is 
small birational, with the exceptional locus $E=f^{-1}P$, 
and crepant; it exists by
\cite[Theorem~0.2]{R83} because we can assume that 
$a(P)=a(Z,0,P)> 1$ (otherwise $a(P)\le1< (n+1)/2$
for $n\ge2$).  
Thus $a(P)=a(Z,0,P)=a=a(X,0,E)$.
We can suppose also that
$X$ is nonsingular and $d\ge 1$ since otherwise, by 
the above $\mathbb Q$-factorial case, $a\le n/2$.
If, in addition, it is extremal then, by
Example~\ref{wish},
the arguments in the proof of our Theorem,
Example~\ref{tor}, and 
the existence of any $D$-flip 
(flop) \cite[Example~3.5.1]{Sh00},
$d\ge (n-1)/2$ (cf. Corollary~\ref{torc}), 
and $a\le n-d\le (n+1)/2$.
The case $a=(n+1)/2$ is possible only when
$d=(n-1)/2$, and 
this is a minimal contraction. 
Since the latter is toric and extremal,
$E$ is a projective nonsingular toric variety
with only ample invariant divisors, that is,
$E=\mathbb P^d$ (the Fano variety of the maximal
index $d+1$), and this is the above contraction.
Suppose now that the relative Picard number of $X/Z$
is $2$. 
Then $X/Z$ has two extremal contractions$/Z$
(it is known in the toric geometry, and 
follows from the Cone Theorem \cite[Theorem~4-2-1]{KMM}
after a perturbation of $B=0$).
They are small and, by the induction on  dimension $n$, 
give $a\le n/2$ if one of them does not contract the
exceptional locus into a point.
Thus, if $a\ge (n+1)/2$, 
we can assume that both are again the above
contractions of $E_1=\mathbb P^d\not=E_2=\mathbb P^d$
with $d=(n-1)/2$, and, by the extremal properties
of both contractions $E_1\cap E_2$ is at
most a point. 
Since $E$ is connected and
all numerically negative
invariant divisors on $E_i/Z$ give $E_i$ in intersection,
there is one of them, say $D_i$, 
which is numerically positive$/Z$ on 
the other $E_{3-i}$. According to the previous
description $D_1+D_2\equiv 0/Z$, and their
intersection $D_1\cap D_2$ gives the induction
in dimension $n$. However, such a contraction
is impossible when $n=3$.
Indeed, for $n=3$, there is no invariant divisor $D$ which
is numerically negative$/Z$ on two curves in $E$, because 
these curves form a connected exceptional sublocus of
$E$ in $D$ and each of its components is a $(-1)$-curve.
In this situation, each $E_i$ is an intersection
of two pairwise distinct invariant divisors.
So, $E=E_1\cup E_2$ because each invariant divisor
passing through a curve of $E$ is negative on
this curve only. 
This gives a contradiction because 
$4$ invariant divisors passing through one of
curves $E_i$ pass through
the intersection point $E_1\cap E_2$.
In particular, we prove that $a\le n/2$ when
the relative Picard number is $2$.
The same holds for higher Picard numbers
by the last case and the induction 
when we consider a contraction
of $2$ dimensional face of the Kleiman-Mori cone
for $X/Z$. 
This completes our proof and explains what are
the minimal toric contractions.\footnote{A.~Borisov
knows a pure combinatorial proof of this fact.} 

So, in dimension $n=3$, Atiyah's flops (2-3) and
their contractions are the only {\em
nonsingular toric minimal\/} ones.

Other elementary examples of minimal semistable flops are
Reid's pagodas and their possible higher 
dimensional generalizations
(cf. La Torre Pendante \cite[8.12]{Ka}).
For $n\ge 5$, the induced flop on $H$ can be different
from the toric induced (or symplectic) one.
All nonsingular minimal flops and
their contractions for $3$-fold are semistable
(that is, their contractions $Z$ have
the cDV type).
They are {\em absolutely\/} extremal, that is,
formally$/\overline{k}$ or analytically when $\overline{k}=\mathbb C$,
if and only if $E=\mathbb P^1$ is irreducible, and
the flops are pagodas -- (1-3) (5) above.
Other minimal $3$-fold flops are their composite
over $\overline{k}$, e.g.,
modifications of Kulikov's model for
K3 surfaces semistable degenerations.
\end{example}

\begin{question}
\label{nonp}
Does Conjecture~\ref{mineq} hold
for nonprojective $X/Z$?
Or at least in the nonsingular case of
Example~\ref{minc}?
It is interesting and nontrivial for $n\ge 4$
because in dimension $3$ each small $X/Z$ is
at least formally or analytically projective.
\end{question}

\begin{question}
\label{minf}
(Cf. Example~\ref{minc}.)
What do minimal contractions of nonsingular $X$
with $n\ge 5$ look like?
Their flops? Absolutely extremal?
Are they still semistable? Their combinatorics?
Is $E$ irreducible and normal ($=\mathbb P^d$)
for absolutely extremal contraction $X/Z$?
\end{question}

\begin{example}[Almost minimal contractions]\label{amc}
Under conditions (BIR) and (WLF),
the next important class having pure $d=a$ is
{\em almost minimal}.

Suppose also that $X$ is nonsingular and $B=0$ then,
for such a contraction, 
$d=a=n/2$, and the dimension is pure, but
either $0$-log or non$0$-log pairs
$(X/Z,0)$ are possible, e.g., Mukai's flop and Kawamata's
flip for $n=4$.
By the theorem it is expected  that $f(E)\le 1$, and
$f(E)$ is a nonsingular curve in $P\in f(E)$ only when
$(X/Z,0)$ is locally over $Z$ a $0$-log pair
with absolutely extremal $X/Z$ near $P$;
then it is expected to be a minimal nonsingular contraction
over the transversal  hyperplane sections
of $f(E)$ through $P$.
The directed or projective flop is fibred.
In particular, this should be a nonsingular
(nonsymplectic even locally over the contraction,
for example, with $\mathbb P^{d-1}$-fibration because
then the lines have $2d-4+1=2d-3=n-3<n-2$
parameters and this contradicts to
Ran's estimation \cite[Lemma~2.3]{bhl})
flop.

The case with the point $f(E)$ is more complicated.
First, suppose that $(X/Z,0)$ is a non$0$-log pair
(this is the {\em minimal\/} case
{\em among\/} the {\em non$0$-contractions\/};
cf. Question~\ref{minn} below).
Then in dimension $4$ the MMP implies that
the flip $X-\to X^+$ exists, and $X^+/Z$ is
nonsingular with a curve ${^+}E=E^+$.
The absolutely extremal components of such flips are
Kawamata's flips (7).
One can hope for a similar picture in
higher dimensions.

Now suppose that $f(E)=pt.$ is a closed point, and
$(X/Z,0)$ is a $0$-log pair.
Such contractions appear as the toric ones induced
in Example~\ref{minc} and also as
the {\em nonsingular\/} birational {\em symplectic\/}
contractions with pure dimension $\le n/2$ for $E$
(then the pure dimension of $E$ is $n/2$ and
as we know $E/pt$).
In these cases, for $n\ge 4$, one can expect the existence
of a nontrivial nonsingular "symmetric"
direct and/or projective flop.
It should have $E^+={^+}E$ of
the pure dimension $d=n/2$ with the same number
of irreducible components as $E$ (the number
exceptional prime divisors over $E$ or
over $E^+$ with the log discrepancy $a=n/2$)
in the formal or analytic case.
For instance, conjecturally the Mukai flop is
the only
nonsingular almost minimal symplectic flop
in dimension $n=4$
(cf. Question~\ref{amin} below).

But there are also singular flops.
For instance, in dimension $4$, an extremal flop
$X-\to X^+/Z$ which transforms $E=\mathbb P^2$
into ${^+}E=E^+=\mathbb P^1$ with a single simple
singularity having the m.l.d. $=2$.

In dimension $2$, the almost minimal contractions
$(X/Z,0)$, which are $0$-log pairs, are
the minimal resolutions of Du Val singularities
$P=f(E)\in Z$.
They are always symplectic, but toric only of
type $\mathbb A_*$, and,
by Example~\ref{minc}, induced toric 
only of type $\mathbb A_1$ (correspond to
the ordinary double singularity). 
\end{example}

\begin{question}
\label{minn}
Is the nonsingular non$0$-log pair case
with $d=n/2$ in higher dimensions similar to
dimension $n=4$?
In particular, is the flip $X^+$ always nonsingular?
$E=\mathbb P^d$ for formally or analytically
extremal contractions?
\end{question}

\begin{question}
\label{amin}
Classify the nonsingular birational almost minimal
$0$-log pairs $(X/Z,0)$ with $f(E)=pt.$ and their
flops.
In particular, such nonsingular symplectic flops.
Is any such flop induced toric 
when $n\ge 4$ (see Example~\ref{minc})?
\end{question}

\begin{example}[More minimal contractions]
A {\em nonidentical\/} contraction
with the m.l.d. $a\le 1$ is never minimal.
For the next (terminal) segment $a\in (1,2]$,
the contraction is minimal only when $d=1$ and
$E$ has only curve components $C$.
Moreover, $X/Z$ is {\em small\/} for $n\ge 3$.
In particular, in dimension $3$,
$E=C$ is a curve.

Such contractions for terminal $3$-folds with $B=0$
appear in the MMP as $0$-log pairs when $a=2$ and
non$0$-log pairs with $a=(m+1)/m$ where $m\ge 2$ is
the index of $K$ in $E=C$.
So, $a\le 3/2$ in the latter case.
In addition, Francia's flip of Example~\ref{frc}
corresponds to (the only one with difficulty $1$)
an extremal
terminal minimal log contraction $(X/Z,0)$
with $a=3/2$ or of the index $2$; it also has 
the maximal length $l=1/2$ among the index $2$
contractions (cf. Example~\ref{benv}).
For $3$-folds with locally complete intersection singularities,
we have only $0$-log pairs with terminal
Gorenstein singularities and their flops that are well-known.

Flops similar to the latter in dimension $n=4$
are still not classified (even the absolutely extremal 
amongst them;
cf. Question~\ref{amin}).
However flips of some (maybe nonminimal)
terminal Gorenstein contractions are explicitly known (8);
they have $d=2$ by Example~\ref{benv}.
\end{example}

\begin{adv}\label{acclt}
An opposite class of {\em maximal\/}
contractions and its application to
the termination of log flips will be given
in one of our future letters.
\end{adv}

\subsection*{List of notation and terminology}

\medskip
\noindent
$a=a(E)=a(X,B,E)$, the m.l.d. of 
$(X,B)$ {\em in\/} the exceptional locus $E$

\medskip
\noindent
$a(Y)=a(X,B,Y)$, the {\em m.l.d.
(minimal log discrepancy)\/} of $(X,B)$ {\em in\/}
subvariety $Y$,
that is, the minimal
log discrepancy $a(X,B,P)$ at
the prime bi-divisors $P$
having the center {\em in\/} $Y$;
the latter means that
$\cent{X}{P}\subseteq Y$;
this assumes that $Y\not=\emptyset$
(otherwise we can put $a(Y)=-\infty$; cf.
Conjecture~\ref{mineq})

\medskip
\noindent
$a(X,B,P)$, for a prime bi-divisor $P$
(prime divisors on some model of $X$ \cite[p.~2668]{Sh96}),
the log discrepancy of $(X,B)$ or $K+B$ 
at $P$ \cite[p.~98]{Sh92}; $P$ is considered
here as its general or scheme point but not as
a subvariety  

\medskip
\noindent
$B$, a Weil $\mathbb R$-divisor on $X$; usually
a boundary (except for Monotonicity), that is,
all its multiplicities $0\le b_i\le 1$, and
$K+B$ is $\mathbb R$-Cartier

\medskip
\noindent
$B^+$, a Weil $\mathbb R$-divisor on $X^+$
such that $f_*^+B^+=f_*B$; usually
a boundary (except for Monotonicity);
for the log flips,
the birational transform of $B$

\medskip 
\noindent
(BIR), the condition on p.~\pageref{bir} which
we assume afterwards, e.g., in Conjecture~\ref{mineq}

\medskip
\noindent
${^+}c$, the {\rm minimal\/} codimension
in $X^+$ of the irreducible components of
the {\rm rational transform\/}
${^+}E$ of $E$ in $X^+/Z$; the codimension is
{\em pure} when all the irreducible
components are of the same codimension

\medskip
\noindent
$\cent{X}{P}$, for a prime bi-divisor $P$
(prime divisors on some model of $X$ \cite[p.~2668]{Sh96}),
its center in $X$ \cite[p.~2669]{Sh96};$P$ is considered
here as a subvariety   

\medskip
\noindent
(CDM), the property on p.~\pageref{cdm} under
additional assumptions in Conjecture~\ref{mineq}

\medskip
\noindent
$\codim P$, for a scheme point $P\in X$,
its codimension $\dim X-\dim P$ in $X$, e.g.,
$\codim P=\dim X$ if and only if $P$ is
a closed point in $X$

\medskip
\noindent
$d=d(X/Z)$, the {\rm minimal\/} dimension
of the irreducible components of
the exceptional locus $E$ of $X/Z$;
this assumes that $E\not=\emptyset$
(otherwise we can put $d=-\infty$; cf.
Conjecture~\ref{mineq}); the dimension is
{\em pure} when all the irreducible
components are of the same dimension

\medskip
\noindent
$D$, a Weil $\mathbb R$-divisor $D$ on $X$

\medskip
\noindent
$D^+$, a semiample$/Z$ $\mathbb R$-Cartier divisor 
on $X^+$ \cite[Definition~2.5]{Sh96} such that
$f_*^+D^+\sim_{\mathbb R} f_*D$  (see Definition 
on p.~\pageref{qflip}); for the log flips,
$D^+$ is the birational transform of $D$
\cite[p.~98]{Sh92}

\medskip
\noindent
$E=E(X/Z)$, the exceptional locus $E$ of $X/Z$, that is,
the union of contractible curves

\medskip
\noindent
$E^+=E(X^+/Z)$, the exceptional locus of $X^+/Z$;
thus $E$ means here to be exceptional;  
in general, ${^+}E$ is quite different from
the rational transform ${^+}E$ 

\medskip
\noindent
${^+}E$, the {\em rational\/} or {\em complete birational
transform\/} of $E$ in $X^+/Z$ (see ${^+}Y$)

\medskip
\noindent
$f:X\to Z$, a birational contraction of normal algebraic
varieties over $k$; sometimes we use such contractions in
the analytic or formal category (cf.
Example~\ref{benv}), and
most of the statements work in them

\medskip
\noindent
$f^+:X^+\to Z$, its qflip or log flip

\medskip
\noindent
$K, K_{X^+}$,
canonical divisors respectively on $X$ and $X^+$
given by the same differential form,
or by the same bi-divisor \cite[Example~1.1.3]{Sh96}

\medskip
\noindent
$l=l(X/Z,B)$, the {\em (log) length\/}
of log pair $(X/Z,B)$; see p.~\pageref{logl}

\medskip
\noindent
LMMP, the log minimal model program and
its conjectures \cite[5.1]{Sh96}

\medskip
\noindent
$n=\dim X$, the dimension of $X$

\medskip
\noindent
(NSN), the property on p.~\pageref{nsn} under
additional assumptions in Conjecture~\ref{mineq}

\medskip
\noindent
(PDM), the property on p.~\pageref{pdm} under
additional assumptions in Conjecture~\ref{mineq}

\medskip
\noindent
$pt.$, a closed point

\medskip 
\noindent
(WLF), the condition on p.~\pageref{wlf} which
we assume afterwards, e.g., in Conjecture~\ref{mineq}

\medskip
\noindent
${^+}Y$, the {\em rational\/} or {\em complete birational
transform\/} of $Y$ in $X^+/Z$; it is defined for
any subvariety $Y\subseteq X$ and any rational map
$g:X-\to Y$ as $g(Y)=\psi\circ\phi^{-1}(Y)$, where
$$
g=\psi\circ\phi^{-1}:X\stackrel{\phi}\gets W
\stackrel{\psi}\to Y
$$
with a birational contraction $\phi$, and 
independent on the decomposition

\medskip
\noindent
$X-\to X^+/Z$, either a $D$-qflip$/Z$, or log 
qflip (see Definition 
on p.~\pageref{qflip}), or a log flip \cite[p.~2684]{Sh96}

\medskip
\noindent
$(X,B)-\to (X^+,B^+)/Z$, in Monotonicity,
a qflip in lf (see Definition 
on p.~\pageref{qflip}) with possibly
nonboundaries $B$ and $B^+$

\medskip
\noindent
$(X/Z,B)$, a log pair which usually satisfies
(BIR) and (WLF) on p.~\pageref{bir}

\medskip
\noindent
a $0$-log pair is a log pair $(X/Z,B)$ such that
$(X,B)$ is log canonical and $K+B\equiv 0/Z$
(cf. \cite[Remark~3.27, (2)]{Sh00})

\medskip
\noindent
$\sim_{\mathbb R}$, the $\mathbb R$-linear
equivalence \cite[Definition~2.5]{Sh96}

\bigskip
\noindent Department of Mathematics,\\ Johns Hopkins University,\\
Baltimore, MD--21218, USA\\ e-mail: shokurov@math.jhu.edu


\begin{thebibliography}{99}
\bibitem{am} Ambro~F., {\em On minimal log discrepancies\/},
Math. Res. Letters {\bf6} (1999), 573--580
\bibitem{AW} Andreatta~M., Wi\'sniewski~J.,
{\em A view on contractions of higher
dimensional varieties\/},
Proc. of Symp. in Pure Math.
{\bf 62.1} (1997), 153--183
\bibitem{at} Atiyah~M. F., {\em On analytic surfaces
with double points\/}, Proc. Roy. Soc. London.
Ser. A {\bf247} (1958), 237--244
\bibitem{Be} Benveniste~X.,
{\em Sur les c\^one des 1-cycles effectifs en
dimension 3\/}, Math. Ann. {\bf 272\/}
(1985), 257--265
\bibitem{BH} Boden~H., Hu~Y.,
{\em Variations of moduli of
parabolic bundles\/}, Math. Ann.
{\bf 301} (1995), 539--559
\bibitem{bhl} Burns~D., Yi Hu, Tie Luo,
{\em HyperK\"ahler manifolds and birational
transformations in dimension 4\/},
preprint.
\bibitem{ChP} Cheltsov~I., Park~J.,
{\em Generalized Eckardt points\/},
submitted for publication (e-print:
Math.AG/0003121).
\bibitem{F} Francia~P.,
{\em Some remarks on minimal models\/}. I,
Compositio Math. {\bf 40} (1980), 301--313
\bibitem{Isk} Iskovskikh~V.A., {\em Birational
rigidity of Fano hypersurfaces in
the framework of Mori theory\/},
Russian Math. Surveys {\bf 56} (2001), 207--291
\bibitem{IP} Iskovskikh~V.A.,
Prokhorov~Yu.G., {\em Fano varieties\/},
Encycl. of Math. Sciences, {\bf 47} (1999),
Springer, Berlin.
\bibitem{Ka} Kachi~Y., {\em Flips from 4-folds
with isolated complete intersection singularities\/},
Amer. J. Math. {\bf 120} (1998), 43--102
\bibitem{K88} Kawamata~Y., {\em Crepant blowing-up of
3-dimensional canonical singularities and
its application to degenerations of
surfaces\/}, Ann. of Math. {\bf 127} (1988), 93--163
\bibitem{K89}  Kawamata~Y., {\em Small contractions
of four dimensional algebraic manifolds\/},
Math. Ann. {\bf 284} (1989), 595--600
\bibitem{KMM} Kawamata~Y., Matsuda~K.,
Matsuki~K., {\em Introduction to
the  Minimal Model Problem\/},
Algebraic Geom., Sendai, 1985, T.~Oda , editor,
Adv. Stud. Pure Math. {\bf 10} (1987),
Kinokuniya Book Store, Tokyo and
North Holland, Amsterdam, 283-360
\bibitem{Ku} Kulikov~Vik.S.,
{\em Degenerations of K3 surfaces and
Enriques surfaces\/}, Math. SSSR
Izvestija {\bf 11\/} (1977),
957--989
\bibitem{M} Markushevich~D., {\em Minimal discrepancy
for a terminal cDV singularity is 1\/},
J. Math. Sci. Univ. Tokyo {\bf3} (1996), 445-456
\bibitem{R83} Reid~M., {\em Decomposition of
toric morphisms\/}, Arithmetic and Geometry II,
M.~Artin and J.~Tate, editors,
Progress in Math. {\bf 36} (1983), Birkh\"auser, 395--418
\bibitem{R} Reid~M., {\em Minimal models of
canonical $3$-folds\/},
Algebraic Varieties and Analytic Varieties,
S.~Iitaka, editor,
Adv. Stud. Pure Math. {\bf 1} (1983),
Kinokuniya Book Store, Tokyo and
North Holland, Amsterdam, 131--180
\bibitem{Sh83} Shokurov~V.V.,
{\em The nonvanishing theorem\/},
Math. USSR Izvestija {\bf 26} (1986),
591--604
\bibitem{Sh88} Shokurov~V.V.,
{\em Problems about Fano varieties\/},
in Birational Geometry of Algebraic Varieties:
Open problems.
The XXIIIrd International Symposium, Division of
Mathematics, The Taniguchi Foundation.
August 22 -- August 27, 1988, 30--32
\bibitem{Sh92} Shokurov~V.V., {\em 3-fold log
flips\/}, Russian Acad. Sci. Izv. Math.
{\bf 40} (1993), 95--202
\bibitem{Sh94} Shokurov~V.V., {\em Semistable
$3$-fold flips\/}, Russian Acad. Sci. Izv. Math.
{\bf 42} (1994), 371--425
\bibitem{Sh96} Shokurov~V.V.,
{\em 3-fold log models\/},
J. Math. Sci.
{\bf 81\/} (1996), 2667--2699
\bibitem {Sh96b} Shokurov~V.V., {\em
Anticanonical boundedness
for curves\/},
appendix to Nikulin~V.V. ``The diagram method for 3-folds and
its application to the K\"ahler cone and
Picard number of Calabi-Yau 3-folds''
in Higher-dimensional complex varieties (Trento,
1994) ed. M. Andreatta, T. Peternell; Berlin: de
Gruyter, 1996, 321--328
\bibitem{sh95} Shokurov~V.V., {\em Complements on surfaces\/},
J. Math. Sci. {\bf102} (2000), 3876--3932
\bibitem{Sh00} Shokurov~V.V.,
{\em Prelimiting flips\/}, preprint, Baltimore-Moscow
(available on
http://www.maths.warwick.ac.uk/~miles/Unpub/Shok/pl.ps), 2001, 235pp.
\bibitem{Th} Thaddeus~M., {\em Geometric invariant theory and
flips\/},
J. Amer. Math. Soc. {\bf9} (1996), 691--723
\bibitem{T} Tsunoda~S., {\em Degenerations of
surfaces\/}, Algebraic Geom.,
Sendai, 1985, T.~Oda , editor,
Adv. Stud. Pure Math. {\bf 10} (1987),
Kinokuniya Book Store, Tokyo and
North Holland, Amsterdam, 755--764
\bibitem{W}  Wi\'sniewski~J.,
{\em On contractions of extremal rays of Fano manifolds\/},
J. reine. angew. Math. {\bf 417\/} (1991), 141--157
\end{thebibliography}
\end{document}